\documentclass[11pt]{article}

\usepackage[english]{babel}
\usepackage{amsmath,amssymb,amsthm,amsfonts}
\usepackage{geometry}
\usepackage{biblatex}
\usepackage{graphicx}
\newtheorem{theorem}{Theorem}[section]
\newtheorem{corollary}[theorem]{Corollary}
\newtheorem{proposition}[theorem]{Proposition}
\newtheorem{remark}[theorem]{Remark}

\theoremstyle{definition}
\newtheorem{definition}{Definition}[section]
\title{Gödel--Dummett and $\mathsf{BD_2}$:\\
Linearity and Depth-Two Branching in Kripke Semantics}
\author{Vicent Navarro Arroyo}
\date{}

\begin{document}

\maketitle

\begin{abstract}
    We study the semantic relationship between Gödel-Dummett logic $\mathsf{GL}$ and bounded-depth-2 logic $\mathsf{BD_2}$, two well-known intermediate logics. While $\mathsf{GL}$ imposes linearity on Kripke frames, $\mathsf{BD_2}$ bounds their depth to two. We prove these logics are incomparable (neither contains the other) through minimal frame conditions. Notably, their combination $\mathsf{GL+BD_2}$ collapses to the logic of one or two world frames, bringing it remarkably close to classical logic. This illustrates how controlling breadth and depth in intuitionistic semantics leads to mutually exclusive structural constraints. Finally, we give a conceptual and philosophical interpretation of the previous results.
This is an extended abstract of work in progress. Comments and suggestions welcome at: this http URL@ub.edu
\end{abstract}

\section{Introduction}

Intermediate logics are those propositional logics that consistently extend intuitionistic logic. They lie between intuitionistic logic and classical logic, since classical propositional logic is the strongest extension of intuitionistic logic, i.e., it proves the most results \cite{sep-logic-intuitionistic}. For example, Gödel--Dummett logic ($\mathsf{GL}$) imposes a global linearity condition on Kripke frames and is axiomatized by the schema
\[(A\to B)\vee(B\to A).\]
Another interesting example is the logic $\mathsf{BD_2}$ which requires that Kripke frames have, at most, branches of length two \cite{Gabbay1986-GABSII}. Moreover, it is axiomatized by the schema
\[A\vee\big(A\to(B\vee\neg B)\big).\]
In this work we will study the semantic relationship between both logics. First, we will see that they are incomparable and we will explain this incomparability through their frame conditions: $\mathsf{GL}$ requires linear frames and $\mathsf{BD_2}$ demands frames with, at most, two branches. Furthermore, we will see that these conditions are minimal in a concrete semantic sense.

\section{Semantic Preliminaries}

In this section we will define the notation and definitions needed in the rest of the work precisely. To study both $\mathsf{GL}$ and $\mathsf{BD_2}$ semantically we need to define intuitionistic propositional logic $\mathsf{IPC}$.

The language we will use is propositional and is generated by a countable set of propositions $\mathsf{Prop}$. The formulas we will consider will be denoted by $\mathsf{Form}$ and are generated, using Backus-Naur notation, as
\[
\mathcal F \ \ := \ \ \top \mid \bot \mid  \mathsf{Prop} \mid (\mathcal F \wedge \mathcal F) \mid (\mathcal F \vee \mathcal F) \mid (\mathcal F \to \mathcal F).
\]
Additionally, we will use the notation $\neg A\colon= A\to\bot$ which designates the usual definition of intuitionistic negation.
\begin{definition}[$\mathsf{IPC}$]\label{definition:IPC}
The axioms and rules of intuitionistic propositional logic ($\mathsf{IPC}$) are standard (see, for example, \cite{sep-logic-intuitionistic}). We use the usual propositional language with connectives $\land$, $\lor$, $\to$, $\bot$, defining $\neg A := A \to \bot$. The only inference rule is \emph{modus ponens}.
\end{definition}

\begin{definition}[$\mathsf{CPC}$]
Classical logic will be denoted by $\mathsf{CPC}$ and arises from $\mathsf{IPC}$ by adding the law of excluded middle $A \vee \neg A$.
\end{definition}
In fact, $\mathsf{CPC}$ is sound and complete with respect to the class of trivial frames (with a single world) \cite{on-the-contingency}.

To study the semantic relationship between $\mathsf{GL}$ and $\mathsf{BD_2}$ we need to define their semantics which will be based on intuitionistic Kripke semantics.

\begin{definition}[Kripke frames and models]
    A Kripke frame $\mathcal{F}$ is a tuple $(W,\leq)$ where $W$ is a non-empty set of worlds and $\leq$ is a partial order relation on $W$. 

    A Kripke model $\mathcal{M}$ is a tuple $(W,\leq, V)$ where $V\colon W\to\mathcal{P}(\mathsf{Prop})$ is a monotone valuation, that is, it satisfies that if $x\leq y$, then $V(x)\subseteq V(y)$.
\end{definition}

Given a Kripke model $\mathcal{M}$, the forcing relation on intuitionistic Kripke models is defined inductively in the standard way.
\begin{align*}
    \mathcal{M}, x\Vdash p &\ \ :\Leftrightarrow \ \ p \,{ \in }\, V(x);\\
    \mathcal{M},x\nVdash \bot;\\
    \mathcal{M}, x\Vdash \varphi \vee \psi&\ \ :\Leftrightarrow \ \ \mathcal{M}, x\Vdash \varphi \mbox{ or }\mathcal{M}, x\Vdash \psi;\\
    \mathcal{M}, x\Vdash \varphi \wedge \psi&\ \ :\Leftrightarrow \ \ \mathcal{M}, x\Vdash \varphi \mbox{ and }\mathcal{M}, x\Vdash \psi;\\
    \mathcal{M}, x\Vdash \varphi \to \psi&\ \ :\Leftrightarrow \ \ \forall \, y{\geq}x\, (\mbox{ if } \mathcal{M}, y\Vdash \varphi \mbox{, then }\mathcal{M}, y\Vdash \psi).
\end{align*}
For simplicity, if the Kripke model is known, we will write $x\Vdash\varphi$ instead of $\mathcal{M},x\Vdash \varphi$.

We say that a formula $\varphi$ is \emph{valid} in a frame $\mathcal{F}$, denoted as $\mathcal{F}\vDash\varphi$, if it is forced at all worlds under
every valuation, that is,
\[\mathcal{F}\vDash\varphi \quad \ \ :\Leftrightarrow\quad \forall x\,\forall V\,(\mathcal{F},V),x\Vdash\varphi.\]
Additionally, recall that the height (or depth) of a frame $\mathcal{(W,\leq)}$ is the supremum of the lengths of all $\leq$-chains over $W$. Conversely, the width of a frame $(W,\leq)$ refers to the number of pairwise incomparable worlds in $W$. 

Finally, we define what a frame condition is since it will help us characterize a logic semantically \cite{Joost_thesis}.
\begin{definition}
    Let $\mathsf{X}$ be a schema of propositional logic. We say that a formula $\mathcal{C}$ in first-order or higher-order predicate logic is a frame condition for $\mathsf{X}$ if
    \[\mathcal{F}\vDash\mathcal{C}\quad\Leftrightarrow\quad \mathcal{F}\vDash \mathsf{X}.\]
\end{definition}
We observe that $\mathcal{F}\vDash \mathcal{C}$ is evaluated in the standard way. In this case $\mathcal{F}$ is understood as a structure in first-order or higher-order predicate logic.

\section{Gödel--Dummett Logic}

We define the logic $\mathsf{GL}$ which extends $\mathsf{IPC}$ and requires linear Kripke frames.

\begin{definition}[$\mathsf{GL}$]
    Gödel--Dummett logic ($\mathsf{GL}$) is the intermediate logic that extends $\mathsf{IPC}$ with the schema
    \[(A\to B)\vee (B\to A).\]
\end{definition}
It is known that $\mathsf{GL}$ is sound and complete with respect to the class of linear Kripke frames \cite{Dummett1959-DUMAPC}.

\begin{theorem}[Frame condition for $\mathsf{GL}$]\label{teorema: fc_GL}
    Given a Kripke frame $\mathcal{F}=(W,\leq)$,
    \[\mathcal{F}\vDash(p\to q)\vee(q\to p) \Leftrightarrow \mathcal{F}\vDash \forall x\,\forall y\,\forall z\,\big((x\leq y\,\wedge\, x\leq z )\to( y\leq z \,\vee \, z\leq y)\big).\]
\end{theorem}
\begin{proof}
    \begin{itemize}
        \item[] 
        \item[$\boxed{\Leftarrow}$] Suppose the frame relation is total. Let us see that for every $x\in W$ if $x\nVdash p\to q$, then $x\Vdash q\to p$.

        If $x\nVdash p\to q$, then there exists some $y\in W$ such that $x\leq y$ and $y\Vdash p$, but $y\nVdash q$. Now suppose that for some arbitrary $z\in W$, $x\leq z$ and $z\Vdash q$. We need to show that $x\Vdash p$. By the frame condition, $y\leq z$ or $z\leq y$, but if $z\leq y$, then (by monotonicity of $\leq$) we would have $y\Vdash q$ which contradicts $y\nVdash q$. Therefore, $y\leq z$ and, again by monotonicity, $z\Vdash p$.
        
        \item[$\boxed{\Rightarrow}$] Suppose that $\mathcal{F}$ does not satisfy the frame condition. Thus, there exist $x,y,z\in W$ such that $x\leq y,\, x\leq z, y\not\leq z$ and $z\not\leq y$. To show that $\mathcal{F}\nvDash (p\to q)\vee(q\to p)$ we must find a valuation such that there is a world in the frame that does not satisfy $(p\to q)\vee(q\to p)$. Consider the following valuation.
        \begin{align*}
            V(p)\colon=\{w\in W\,\colon\, w\leq y\}, & &V(q)\colon=\{w\in W\,\colon\, z\leq w\}.
        \end{align*}
        It is easy to see that $x\nVdash p\to q$ and $x\nVdash q\to p$.
    \end{itemize}
\end{proof}
Observe that this frame condition is minimal and sufficient. Moreover, it has been enough to find a single branching in the frame to demonstrate the failure of the schema $(A\to B)\vee(B\to A)$.

\section{The logic $\mathsf{BD_2}$}

We define the logic $\mathsf{BD_2}$ which extends $\mathsf{IPC}$ requiring frames of depth at most $2$.

\begin{definition}[$\mathsf{BD_2}$]
    Bounded depth $2$ logic ($\mathsf{BD_2}$) is the intermediate logic that extends $\mathsf{IPC}$ with the schema
    \[A\vee(A\to (B\vee \neg B)).\]
\end{definition}
It is well known that $\mathsf{BD_2}$ is sound and complete with respect to the class of frames of depth at most two \cite{Gabbay1986-GABSII}.

\begin{theorem}[Frame condition for $\mathsf{BD_2}$]\label{teorema: fc_BD2}
    Given a Kripke frame $\mathcal{F}=(W,\leq)$,
    \[\mathcal{F}\vDash p\vee\big(p\to(q\vee\neg q)\big) \Leftrightarrow \mathcal{F}\vDash \forall x\,\forall y\,\forall z\,\big((x\leq y\,\wedge\, x\leq z )\to\big( y\leq z \,\to \, (y=x \,\vee z=x)\big)\Big).\]
\end{theorem}
\begin{proof}
    \begin{itemize}
        \item[]
        \item[$\boxed{\Leftarrow}$] Suppose the relation allows, at most, two branches. Let us see that for every $x\in W$ if $x\nVdash p$, then $x\Vdash p\to(q\vee\neg q)$.

        Suppose $x\nVdash p$ and let $y\in W$ be arbitrary such that $x\leq y$ and $y\Vdash p$. We need to show that $y\Vdash q \vee \neg q$. Therefore, suppose that $y\nVdash \neg q$ and let us see that $y\Vdash q$. If $y\nVdash \neg q$, then there exists some $z\in W$ such that $y\leq z$ and $z\Vdash q$. By the frame condition, either $y=x$ or $z=x$. However, if $y=x$, then $x\Vdash p$ and $x\nVdash p$, which is contradictory. Therefore, $z=x$ and thus $x\Vdash q$. Finally, by the monotonicity of $\leq$, we conclude that $y\Vdash q$.
        
        \item[$\boxed{\Rightarrow}$] Suppose that $\mathcal{F}$ does not satisfy the frame condition. Thus, there exist $x,y,z\in W$ such that $x < y < z$. To show that $\mathcal{F}\nvDash p\vee\big(p\to (q\vee\neg q)\big)$ we must find a valuation such that there is a world in the frame that does not satisfy $p\vee\big(p\to(q\vee\neg q)\big)$. Consider the following valuation.
        \begin{align*}
            V(p)=\{w\in W\,\colon\, y\leq w\}, &&V(q)=\{w\in W\,\colon\, z\leq w\}.
        \end{align*}
        It is easy to see that $x\nVdash p$ and $x\nVdash p\to(q\vee\neg q)$.
    \end{itemize}
\end{proof}
Note that the frame condition is minimal and sufficient. Moreover, a chain of length $3$ already invalidates the axiom $A\vee\big(A\to(B\vee\neg B)\big)$. That is why this logic is called \textit{bounded depth $2$}.

\section{Semantic separation and incomparability}

Let us see that the logics $\mathsf{GL}$ and $\mathsf{BD_2}$ are incomparable, that is, there exist formulas that are derivable in $\mathsf{GL}$ but not in $\mathsf{BD_2}$ and vice versa. The proof we offer is purely semantic and will depend only on the frame conditions.

\begin{figure}[h!]
    \centering
    \includegraphics[scale=1]{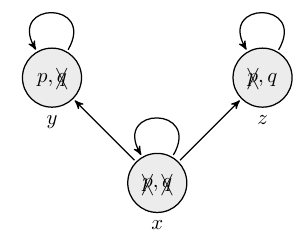}
    \caption{Frame that satisfies $\mathsf{GL}$ but not $\mathsf{BD_2}$. Crossed out propositional letters indicate they are not forced at that world.}
    \label{fig:GL_noBD2}
\end{figure}

\begin{theorem}\label{GL_notBD2}
    $\mathsf{GL}\not\subseteq \mathsf{BD_2}$.
\end{theorem}
\begin{proof}
    We know that $\mathsf{GL}$ and $\mathsf{BD_2}$ are sound and complete with respect to the classes of frames of depth at most two and linear frames, respectively \cite{Dummett1959-DUMAPC, Gabbay1986-GABSII}. Therefore, it will suffice to consider a Kripke frame that satisfies the frame condition of $\mathsf{GL}$ but not that of $\mathsf{BD_2}$ (Theorems \ref{teorema: fc_GL} and \ref{teorema: fc_BD2}). In other words, it suffices to consider a linear frame since the order is total and there are no incomparable worlds.

    Let $\mathcal{F}\colon=(W,\leq)$ such that $W=\{x,y,z\}$ and $x<y<z$. See Figure \ref{fig:GL_noBD2} for a representation of this frame.
    \end{proof} 
   
\begin{figure}[h!]
    \centering
    \includegraphics[scale=1]{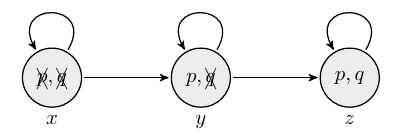}
    \caption{Frame that satisfies $\mathsf{BD_2}$ but not $\mathsf{GL}$. Crossed out propositional letters indicate they are not forced at that world.}
    \label{fig:BD2_noGL}
\end{figure}

\begin{theorem}
    $\mathsf{BD_2}\not\subseteq \mathsf{GL}$.
\end{theorem}
\begin{proof}
We proceed in a similar way to the proof of Theorem \ref{GL_notBD2}. We consider a frame of bounded depth and branched with two branches since there are two incomparable worlds and, thus, the order is not total.

    Let $\mathcal{F}\colon=(W,\leq)$ such that $W=\{x,y,z\}$, $x<y$ and $x<z$. See Figure \ref{fig:BD2_noGL} for a representation of this frame.
\end{proof}

\begin{proposition}[Non-empty intersection]\label{prop: intersec_no_vacia}
    $\mathsf{GL} \cap \mathsf{BD_2} \neq \varnothing$.
\end{proposition}
\begin{proof}
    It suffices to observe that both $\mathsf{GL}$ and $\mathsf{BD_2}$ are consistent extensions of $\mathsf{IPC}$, so every theorem of $\mathsf{IPC}$ belongs to their intersection. Moreover, frames of one or two worlds satisfy both linearity (condition of $\mathsf{GL}$) and depth $\leq 2$ (condition of $\mathsf{BD_2}$), and in them formulas like $\neg\neg(p \vee \neg p)$ are valid which are not theorems of $\mathsf{IPC}$.
\end{proof}

\begin{corollary}[Extension of $\mathsf{GL}$ and $\mathsf{BD_2}$] The logic $\mathsf{GL}+\mathsf{BD_2}$, which extends $\mathsf{IPC}$ with the schemas of $\mathsf{GL}$ and $\mathsf{BD_2}$, is sound and complete with respect to the class of frames of one or two worlds.
\end{corollary}
\begin{proof}
 $\mathsf{GL}+\mathsf{BD_2}$ is the logic of the intersection of the frame classes of $\mathsf{GL}$ and $\mathsf{BD_2}$ and soundness and completeness are preserved under intersection of frame classes.
\end{proof}
\begin{remark}
The logic $\mathsf{GL}+\mathsf{BD_2}$, being that of frames of size $\leq 2$, is very close to $\mathsf{CPC}$. In fact, while $\mathsf{CPC}$ is exactly the logic of the trivial frame, $\mathsf{GL}+\mathsf{BD_2}$ additionally admits the linear frame of two worlds $w_1 < w_2$. This latter frame already validates many classical principles. The progression
\[
\mathsf{IPC} \subsetneq \mathsf{GL},\mathsf{BD_2} \subsetneq \mathsf{GL}+\mathsf{BD_2} \approx \mathsf{CPC}
\]
illustrates how simultaneously imposing orthogonal restrictions on the Kripke structure inevitably leads to classical collapse.
\end{remark}
In summary, $\mathsf{GL}$ and $\mathsf{BD_2}$ are incomparable but not disjoint, and their axiomatic union collapses semantically to frames of size $\leq 2$. This tension between linearity and bounded depth illustrates how apparently natural structural restrictions can be mutually exclusive in the intuitionistic framework.

\section{Conceptual interpretation}

We can give a conceptual interpretation, from an intuitionistic standpoint, of the previous results. To do this we must first understand how to interpret the Kripke frames of $\mathsf{IPC}$. According to the usual interpretation \cite{Troelstra1988-TROCIM}, Kripke models represent states of knowledge that branch into other future states, giving a temporal sense to the partial order relation of the model. The conditional $A\to B$ means that in every future knowledge state I will not be able to know $A$ without knowing $B$. Thus, negating $A$ in a state ($\neg A$) is understood as an exclusion of $A$ in the accessible future states. Moreover, monotonicity ensures that knowledge in a state persists in future states, that is, what is proved remains set in stone.

Linear Kripke frames characterize $\mathsf{GL}$ and, therefore, suggest that in this logic the growth of information towards the future is linear; there cannot exist alternative information states. In fact, the schema $(A\to B)\vee(B\to A)$ indicates that all formulas are ordered by intuitionistic implication. On the other hand, in $\mathsf{BD_2}$, its characterization by depth-bounded frames of two indicates that the justification of knowledge is limited. Its schema $A\vee \big(A\to (B\vee \neg B)\big)$ can be read as: if $A$ is not yet proved, then, when $A$ becomes true, it forces a classical decision about $B$.

Metaphorically, one can understand the philosophical view of $\mathsf{GL}$ on knowledge as that of the optimistic view of refinement of scientific knowledge, since it assumes a single potentially infinite course of refinement that seeks an ideal of wisdom. On the other hand, the view of knowledge of $\mathsf{BD_2}$ would be democratic (with multiple voices), but pragmatic since one must stop after a review. It can be compared with deliberative knowledge in an agora.

Both views represent orthogonal ways of regulating knowledge in intuitionistic logic without being as strict as classical logic. $\mathsf{GL}$ regulates in width (only one voice), while $\mathsf{BD_2}$ regulates in height (does not deepen). In fact, they are not exclusive views since we have already shown that $\mathsf{GL}\cap\mathsf{BD_2}\neq\varnothing$ (Proposition \ref{prop: intersec_no_vacia}), but the combination of both strategies results in the extension $\mathsf{GL}+\mathsf{BD_2}$ which is too strong; it collapses into a logic very close to classical logic since only trivial frames or frames of two worlds are admitted.

\section{Conclusion and future lines}

In this work we have demonstrated the incomparability of $\mathsf{GL}$ and $\mathsf{BD_2}$ in exclusively semantic terms and by using their frame conditions. These frame conditions have been demonstrated explicitly and carefully, as well as illustrated. Moreover, we have revealed how the combination of orthogonal structural restrictions (linearity and bounded depth) inevitably leads to a semantic collapse towards classical logic. Finally, we have offered a conceptual interpretation of the obtained results from the coordinates of intuitionistic logic.

Finally, as future research ideas, one could study the incomparability of other intermediate logics like $\mathsf{BD_n}$, the logics characterized by frames of depth at most $n$. One could also investigate the semantic relations from an algebraic point of view. Additionally, one could investigate the consequences of mixing the logics $\mathsf{GL}$ and $\mathsf{BD_2}$, as is done in \cite{on-the-contingency} with classical and intuitionistic propositional logics.
\newpage

\section{Appendix}

\printbibliography
\end{document}